\documentclass[reqno,11pt]{amsproc}

\usepackage{hyperref}
\usepackage{url}
\newtheorem{theo}{Theorem}
\newtheorem{rem}{Remark}

\newtheorem{lem}{Lemma}

\newtheorem{hyp}{Hypothesis}

\newcommand\eps\varepsilon
\newcommand\ph\varphi
\newcommand\kap\Lambda

\usepackage{enumerate}
\usepackage{graphicx}
\usepackage{float}
\usepackage{yfonts}
\usepackage{mathrsfs}
\usepackage{amsthm}
\usepackage{amsmath}
\usepackage{amscd}
\usepackage{xcolor}

\begin{document}

\title[On Zeros of Homotopic Mappings]
{On Zeros of Homotopic Mappings, Fixed Point Problems, and Inverse Functions}

\author[Oleg Zubelevich]{Oleg Zubelevich\\Faculty of Mechanics and Mathematics,\\
Lomonosov Moscow State University,\\ 1 Leninskie Gory, Moscow 119991, Russia}
\email{oezubel@gmail.com}

\date{}

\begin{abstract}This article examines a family of smooth mappings between Banach spaces and establishes conditions for the existence of their zeros. Applications to fixed-point problems and the Implicit Function Theorem are also discussed.
\end{abstract}

\maketitle
\numberwithin{equation}{section}
\newtheorem{theorem}{Theorem}[section]
\newtheorem{lemma}[theorem]{Lemma}
\newtheorem{definition}{Definition}[section]
\section{Introduction}
The zeros of mappings between Banach spaces have been widely studied from the viewpoint of the continuation method in many articles: \cite{ang}, \cite{[1]}, \cite{[2]}, \cite{[3]}, \cite{[4]}, \cite{[5]}.

An exposition of the continuation method and its
applications can be found in \cite{[22]},  \cite{[23]}, \cite{[24]}.

The continuation method is primarily applied in two contexts: the implicit and inverse function theorems, and fixed-point theory. Some modern works discussing these topics from an applied perspective include \cite{gar} and \cite{man}.

In this article, we develop a continuation method based on the extension of solutions to the associated initial value problem for an ODE. Although the underlying concept has likely always been intuitive, one of the first formal publications to address it appears to be \cite{Yak}. That work, however, is limited to finite-dimensional equations.

From the perspective of differential equations theory, our ability to construct a continuation method is determined by the existence theorems available for the Cauchy problem in various function spaces.

Against this background, existence theorems provided in the following works may be of use:
\cite{safonov}, \cite{mil}, \cite{zy}.

\section{The Main Theorems}
\subsection{The case of general domain}
Let $(X,\|\cdot\|_X),\quad (Y,\|\cdot \|_Y)$ be Banach spaces, and let  $D\subset X$ be an open, non-empty subset. Its boundary is denoted by
$ \partial D=\overline D\backslash D.$

The case $D=X$ and $\partial D = \emptyset$ is possible.

Let
$$B_r(z):=\{x\in X\mid \|x-z\|_X<r\}$$ stand for the  open ball in $X$ centered at $z$ with radius $r$.
Similarly, we denote by
$B_r^Y(\tilde z)$ the  open ball in $Y$ centered at $\tilde z$ with radius $r$.

We use $\mathscr L(X,Y),\quad \mathscr{L}(X):=\mathscr L(X,X)$ to denote the space of bounded linear operators $$A:X\to Y$$ equipped with the uniform operator topology, where the norm is defined as:
$$\|A\|_{\mathscr L(X,Y)}=\sup_{\|x\|_X\le 1}\|Ax\|_Y.$$

Introduce the notation $ I=[0,1)$ and
consider a continuous mapping
$$F=F(t,x),\quad F:\overline I\times\overline D\to Y.$$
\begin{hyp}\label{00}
The mapping $F$ is uniformly continuous in $t$. This means that
for any $\eps>0$, there exists $\delta>0$ such that, if $|t'-t''|<\delta $ then
$$\sup_{x\in D}\|F(t',x)-F(t'',x)\|_Y<\eps,\quad  t',t''\in\overline I.$$
\end{hyp}

The  partial Fr\'echet derivatives
$$F_t(t,x)\in\mathscr L(\mathbb{R},Y),\quad F_x(t,x)\in\mathscr L(X,Y)$$
are defined at each point $(t,x)\in I\times D$.

In what follows, let $E_Y = \mathrm{id}_Y$ and $E_X = \mathrm{id}_X$.

\begin{hyp}\label{h4}
For any $(t,x)\in I\times D$, there exists a right inverse
$$S(t,x)\in \mathscr L(Y,X),\quad F_x(t,x)S(t,x)=E_Y$$
such that the function
$$a:I\times D\to X,\quad a(t,x)=-S(t,x)F_t(t,x)$$ is continuous and locally Lipschitz in $x$ (see Section \ref{cdwwsd}). \end{hyp}
In the simplest case,  $S=F_x^{-1}\in\mathscr{L}(Y,X).$ If $X$ and $Y$ are real Hilbert spaces and the operator $F_x(t,x)$ is onto  for all $(t,x)\in I\times D$,
then the right inverse can be chosen as $S=F_x^*(F_xF_x^*)^{-1}$, where $F_x^*:Y\to X$ is the adjoint operator. A general discussion on the existence of a right inverse operator can be found in \cite{brezis}. Degenerate operators are studied in \cite{AR}.
\begin{hyp}\label{h}
The growth of the function $a$ on any closed set is at most linear.
Specifically, for any closed  set $P\subset D$ there exists a constant $K_P$ such that
$$\|a(t,x)\|_X\le K_P(\|x\|_X+1),\quad (t,x)\in I\times P.$$
\end{hyp}

\begin{theo}\label{q12e}
Under  hypotheses \ref{00}, \ref{h4}, and \ref{h}, suppose that $F(0, x_0)=0$ for some $x_0\in D$. Then,  there are two possibilities:
\begin{enumerate}
  \item \label{step:1}
there exists a number $\tau\in (0,1]$ and a sequence $\{x_k\}\subset\partial D$ such that
$F(\tau,x_k)\to 0$ in $Y$;
and
\item \label{step:2}
 there exists an element $x_1\in D$ such that $F(1,x_1)=0.$\end{enumerate}
\end{theo}

Note that possibility (\ref{step:1}) can be formulated alternatively as:$$0 \in \overline{F(\tau, \partial D)}.$$In this form, it is quite well known in finite-dimensional topological degree theory \cite{nir}.

\subsection{The case of coercive $F$ and an unbounded domain}
Theorem \ref{q12e} holds in the general case, regardless of whether the domain $D$ is bounded or not. If $D$ is bounded, Hypothesis \ref{h} reduces to the assumption that the mapping $a$ is bounded on closed bounded sets. However, we can retain this hypothesis in the case of an unbounded $D$ if we apply the coercivity condition to $F$.

Assume that the domain $D$ is unbounded.
\begin{hyp}\label{u1}
The mapping $F$ is coercive for each $t$. Specifically,  for any sequence $\{x_n\}\subset D$ such that
$\|x_n\|_X\to\infty$, it follows that
$$\|F(t,x_n)\|_Y\to \infty,\quad t\in\overline I.$$
\end{hyp}

\begin{hyp}\label{u2}
The function $a$ is bounded on any closed bounded set $P\subset D$.
Specifically, for any such set $P\subset D$, we have
$$\sup_{(t,x)\in I\times P}\|a(t,x)\|_X<\infty.$$
\end{hyp}
\begin{theo}\label{xfb6t}
Under the hypotheses \ref{00}, \ref{h4}, \ref{u1}, and \ref{u2}, if $F(0, x_0)=0$ for some $x_0 \in D$, then there are two possibilities:

\begin{enumerate}
  \item \label{p:1}
  there exist a number $\tau \in (0,1]$ and a sequence $\{x_k\} \subset \partial D$ such that $F(\tau,x_k) \to 0$ in $Y$;

  \item \label{p:2}
  there exists an element $x_1 \in D$ such that $F(1,x_1) = 0$.
\end{enumerate}
\end{theo}

\section{Applications}\subsection{ A Fixed Point Problem}
We now apply the previous result to a fixed point problem.

 Assume that \begin{equation}\label{xfgbgm}0 \in D.\end{equation}

 Let $f: \overline{D} \to X$ be a bounded continuous mapping that is Fréchet differentiable in $D$, and suppose the map $x \mapsto f'(x)$ is locally Lipschitz on $D$ with values in $\mathscr{L}(X)$.

Under these assumptions the following theorem holds.
\begin{theo}\label{zdsgfv5--}
Assume that for any closed and bounded set $P \subset D$, we have
 \begin{equation}\label{wwxbf700}\sup_{(t,x)\in I\times P}\big\|\big(E_X-tf'(x)\big)^{-1}\big\|_{\mathscr{L}(X)}<\infty,\end{equation}
 and that $f(\partial D)$ is relatively compact.
Furthermore, suppose that for any $x\in\partial D,$ at least one of the following conditions is satisfied:
\begin{enumerate}
\item \label{s:1}
$\|f(x)\|_X\le\|x\|_X;$
\item\label{s:2}$\|f(x)\|_X\le\|x-f(x)\|_X;$
\item\label{s:3}
$\|f(x)\|_X\le\sqrt{\|x\|_X^2+\|x-f(x)\|_X^2}$;
\item\label{s:4}$\|f(x)\|_X\le
\max\{\|x-f(x)\|_X,\|x\|_X\}$.
\end{enumerate}Then $f$ has a fixed point in $\overline D$.
\end{theo}
Note that the inequality (\ref{wwxbf700}) can be replaced with a coarser one:
$$\sup_{x\in P}\|f'(x)\|_{\mathscr{L}(X)} < 1.$$
Observe that if the domain $D$ is not convex, this inequality does not, in general, entail that the mapping $f$ is nonexpansive.

The list of inequalities (\ref{s:1})–(\ref{s:4}) is taken from the book \cite{aga}. A similar list can also be found in \cite{gran}.

If the mapping $f: \overline{D} \to X$ is only continuous, without hypotheses regarding its smoothness, and one of the inequalities  (\ref{s:1})–(\ref{s:4}) for boundary points holds, then the following three types of fixed-point theorems are known.

1. Assume that the mapping $f$ is contractive and the set $f(\overline D)$ is bounded. Then $f$ has a fixed point.

2. If $f$ is a compact map then it has a fixed point (Schauder  type theorem).

3. Assume that $X$ is a Hilbert space, or more generally a uniformly convex Banach space, $D$ is convex and bounded, and $f$ is nonexpansive. Then $f$ has a fixed point (Kirk type theorem).

\subsection*{Proof of Theorem \ref{zdsgfv5--}}
Let us apply Theorem \ref{xfb6t} to the mapping
$$F(t,x)=x-tf(x).$$If the domain $D$ is bounded, then the same argument
applies by Theorem \ref{q12e}.

It is sufficient to handle the  possibility (\ref{p:1}): there exists a sequence
$\{x_k\} \subset\partial D$ and a number $\tau\in(0,1]$
such that$$x_k - \tau f(x_k) \to 0.$$
The sequence $\{f(x_k)\}$ contains a convergent subsequence, say $$f(x'_k)\to \hat x,\quad \{x'_k\}\subset\{x_k\}.$$ It follows that $x'_k\to\tau \hat x\in\partial D$ and
\begin{equation}\label{zsdg66pp}
f(\xi)=\frac{1}{\tau}\xi,\quad \xi=\tau\hat x.\end{equation}
Since $\xi \in \partial D$, by (\ref{xfgbgm}) we have $\xi \ne 0$.

The case $\tau < 1$ is inconsistent with conditions (\ref{s:1})--(\ref{s:4}). Indeed, consider condition (\ref{s:4}), for instance. It then follows from (\ref{zsdg66pp}) that
$$\frac{1}{\tau}\|\xi\|_X\le \max\Big\{\Big\|\xi-\frac{1}{\tau}\xi\Big\|_X,\|\xi\|_X\Big\}
=\max\Big\{\frac{1}{\tau}-1,1\Big\}\|\xi\|_X.$$
This inequality is impossible for $\tau\in(0,1)$.

In the case $\tau = 1$, we obviously obtain $f(\hat{x}) = \hat{x}$. This completes the proof of Theorem \ref{zdsgfv5--}.

\subsection{An Implicit Function Theorem}

Inessential positive constants are denoted by the same letter $c$.

Assume that $f: X \to Y$ is a Fr\'echet differentiable mapping such that the map $x \mapsto (f'(x))^{-1} \in \mathscr{L}(Y, X)$ is locally Lipschitz.
\begin{theo}\label{xdfg..}
Suppose that either$$\sup_{x\in X} \big\|(f'(x))^{-1}\big\|_{\mathscr{L}(Y, X)} \le c(\|x\|_X + 1),$$or the mapping $f$ is coercive and, for any closed and bounded set $P \subset X$,$$\sup_{x\in P} \big\|(f'(x))^{-1}\big\|_{\mathscr{L}(Y, X)} < \infty.$$
Then $f$ has a continuous right inverse $\psi: Y \to X$ such that $f \circ \psi = E_Y$.
\end{theo}
The fact that the mapping $f$ is onto  follows directly from Theorems \ref{q12e} and \ref{xfb6t} by setting
$$F(t,x)=f(x)-(1-t)f(0)-ty,$$where $y\in Y$ is an arbitrary point.

On the other hand, it follows from the proofs of these theorems that the point$$x_1=x_1(y), \quad F(1,x_1(y))=0,$$is obtained via the following procedure.

Consider the initial value problem$$\dot x = (f'(x))^{-1} (y - f(0)), \quad x(0) = 0.$$This problem has a unique solution $x = x(t, y),\quad 0\le t\le 1$. We then set $\psi(y) =x_1(y)= x(1, y)$.
The continuity of the mapping $\psi$ follows from
standard ODE theorems on the continuous dependence of solutions on parameters.

In this regard, recall the following result.
\begin{theo}[\cite{re}]\label{sdggfg4}
Let $G:X\to X$ be compact and continuously differentiable map. Set $\ph=E_X-G$ and suppose that $\ph'$ has a bounded inverse $(\ph'(x))^{-1}\in\mathscr{L}(X)$ such that
$$\|(\ph'(x))^{-1}\|_{\mathscr{L}(X)}\le c(\|x\|_X+1),\quad x\in X.$$
Then, $\ph$ is a homeomorphism, and hence a diffeomorphism, from $X$ to itself.
\end{theo}

\section{Proof of the Theorems  }
\subsection{Proof of Theorem \ref{q12e}}
Suppose that  the possibility \ref{step:1} does not hold, i.e.,
\begin{equation}\label{awee}0\not\in\overline {F(t,\partial D)}\end{equation}
for all $t\in(0,1]$,
and let us show that this implies that the possibility \ref{step:2} is fulfilled.

\begin{lem}\label{xfhb6ooo} For each $t\in  (0,1]$ there exists a positive constant $r$ such that
$$\overline {F^{-1}(t,B^Y_r(0))}\subset D.$$\end{lem}
\subsubsection*{Proof of Lemma \ref{xfhb6ooo}}
By  (\ref{awee}) we can choose $r$ such that
\begin{equation}\label{xcvnn;;}\overline {B^Y_r(0)}\cap \overline{ F(t,\partial D)}=\emptyset.\end{equation}
To prove the lemma, we shall show that $\overline {F^{-1}(t,B^Y_r(0))}\cap \partial D=\emptyset.$

Indeed, assume the converse:
there exists a point $\tilde x\in \partial D$ and a sequence $\{x_n\}\subset{F^{-1}(t,B^Y_r(0))}$ such that $x_n\to\tilde x$.

Thus, we have $F(t,x_n)\in B^Y_r(0)$. Since $F$ is a continuous mapping, it follows that
$$F(t,\tilde x)\in \overline{B^Y_r(0)}.$$ On the other hand,
$$F(t,\tilde x)\in \overline{ F(t,\partial D)}.$$
This contradicts formula (\ref{xcvnn;;}) and the lemma is proved.

Consider an initial value problem
\begin{equation}\label{zsdfv132}
\dot x=a(t,x),\quad x\mid_{t=0}=x_0.\end{equation}
According to Hypothesis \ref{h4} and the Cauchy existence theorem (see Section \ref{cdwwsd} for details), this problem has a unique solution $x(\cdot)\in C^1([0,\tau),X),$  where $\tau\in(0,1]$.

It is easy to see that the function $F(t,x)$ is a first integral of the equation (\ref{zsdfv132}). That is,
\begin{equation}\label{xfbhpp}F(t,x(t))=0\end{equation} for all admissible $t$. Indeed:
\begin{align}\frac{d}{dt}F(t,x(t))=F_t(t,x(t))&+F_x(t,x(t))\dot x\nonumber\\
&=F_t(t,x(t))+F_x(t,x(t))a(t,x(t))=0.\nonumber
\end{align}
On the other hand, from the conditions of the theorem we have
$$F(0,x(0))=0.$$

First assume that $\tau=1$. In accordance with Lemma \ref{xfhb6ooo} choose $r>0$ small enough such that
$$\overline {F^{-1}(1,B^Y_r(0))}\subset D.$$
Observe also that $F^{-1}(1,B^Y_{r}(0))$ is an open set as the preimage of an open set.

By Hypothesis \ref{00} and equality (\ref{xfbhpp}), there is a constant $\delta>0$ such that for each $t\in[1-\delta,1)$,
one has
$$\|F(1,x(t))\|_Y=\|F(t,x(t))-F(1,x(t))\|_Y<r.$$
This implies that
$$x(t)\in  {F^{-1}(1,B^Y_{r}(0))},\quad t\in[1-\delta,1).$$ It follows from Hypotheses \ref{h}, \ref{h4} and  Remark \ref{xdvfcrt00} (see Section \ref{cdwwsd})  that the limit
$$\lim_{t\to 1} x(t)=x_1\in \overline {F^{-1}(1,B^Y_r(0))}$$ exists. Here, the set $F^{-1}(1,B^Y_{r}(0))$ plays the role of $Q$.

It remains to pass to the limit in  equality (\ref{xfbhpp}) to obtain:
$$\lim_{t\to 1}F(t,x(t))=F(1,x_1)=0.$$
This completes the proof of the theorem for the case $\tau=1$.

Suppose that $\tau<1$, and the solution cannot be extended to the right beyond $\tau$. Let us show that this leads to a contradiction.

The further argument repeats the previous one, up to  details.

Indeed, take $r>0$ such that $\overline {F^{-1}(\tau,B^Y_r(0))}\subset D$ and pick $\delta>0$ such that
$$x(t)\in F^{-1}(\tau,B^Y_r(0)),\quad t\in [\tau-\delta,\tau).$$
By Theorem \ref{xdfgfftQ}, the limit
$$\lim_{t\to\tau}x(t)=x_-$$ exists and $x_-\in\partial {F^{-1}(\tau,B^Y_r(0))}\subset D.$

Now, take $r'>0$ such that $\overline {B_{r'}(x_-)}\subset D$. By Theorem \ref{xdfgfftQ} with $Q=B_{r'}(x_-)$, we can extend the solution $x(t)$ to the right beyond $\tau$. This yields a contradiction.

Theorem \ref{q12e} is thus proved.
\subsection{Proof of Theorem \ref{xfb6t}}
The proof of Theorem \ref{xfb6t} essentially repeats the argument from the proof of Theorem \ref{q12e}.

Indeed,  Hypothesis \ref{u1} implies that for any $t\in \overline I$ and  any bounded set $B\subset Y$,
the preimage $F^{-1}(t,B)$ is bounded. Thus, the sets $F^{-1}(1,B^Y_{r}(0))$ and $F^{-1}(\tau,B^Y_{r}(0))$ from the proof of Theorem \ref{q12e} are bounded, and Hypothesis \ref{u2} allows one to employ the results of Section \ref{cdwwsd} in the same way as above.

Theorem \ref{xfb6t} is proved.

\section{Appendix: Extensibility of Solutions to ODE }\label{cdwwsd}
For the sake of completeness, we recall a fundamental theorem from the theory of ordinary differential equations. Various versions of this theorem can be found in standard textbooks \cite{hart}.

Let $Q$ be a nonempty  open subset of a Banach space $X$.

Consider the ordinary differential equation:
\begin{equation}\label{5y}
\dot x=b(t,x).\end{equation}
Assume the mapping $b:[0,T)\times Q\to X$ satisfies the conditions of the Cauchy existence and uniqueness theorem.

Recall these conditions: $b\in C\big([0,T)\times Q, X\big)$ and the mapping $b$ is  locally Lipschitz with respect to $x$. That is, for any point $(\hat t,z)\in [0,T)\times Q$ there exist a ball
$$B_\eps(z)\subset Q,$$ and positive  constants $c,\sigma$ such that
$$\|b(t,x'')-b(t,x')\|_X\le c\|x'-x''\|_X,\quad t\in[\hat t,\hat t+\sigma),\quad x',x''\in B_\eps(z).$$

These conditions guarantee that for any initial value
$$x(t_0)=x_0\in Q,\quad t_0\in[0,T)$$the ordinary differential equation (\ref{5y}) has a unique solution $x(t)$ defined on an interval $t\in[t_0,t_1)$ where   $t_1-t_0>0$ is sufficiently small \cite{swarz}.

\begin{theo}\label{xdfgfftQ}
Assume that there exists a positive constant $K_Q$ such that for any
$(t,x)\in [0,T)\times Q$ we have
$$\|b(t,x)\|_X\le K_Q(\|x\|_X+1)$$ and let
$x(t)$ be a solution to the equation (\ref{5y}) such that $x(t)$ is
defined for $t\in[0,t^*)$ with $ t^*<T$.

Then the limit \begin{equation}\label{zcxv1}\lim_{t\to t^*-}x(t)=x_-\end{equation} exists and there are two possibilities:
\begin{enumerate}
\item $x_-\in Q$, in which case the solution $x(t)$ can  be extended to an interval $[0,\tau)$ with $\tau\in(t^*,T)$;

\item $x_-\in\partial Q$.\end{enumerate}
\end{theo}

\subsubsection*{Proof of  theorem \ref{xdfgfftQ}} We write equation  (\ref{5y}) in the integral form:
$$x(t)=x(0)+\int_0^tb(s,x(s))ds.$$
Thus, given the estimate
$$\|x(t)\|_X\le \|x(0)\|_X+ K_Q\int_0^t(\|x(s)\|_X+1)ds,$$ from  Gronwall's inequality
we obtain:
$$\|x(t)\|_X\le (\|x(0)\|_X+1)e^{K_Qt}-1.$$
This inequality demonstrates that $x(t)$ is bounded on $[0,t^*)$. Thus the function $b(t,x(t))$ is also bounded:
$$\|b(t,x(t))\|_X\le M,\quad t\in [0,t^*).$$
Then
 for $t',t''\in [0,t^*)$ we have the estimate:
$$\|x(t'')-x(t')\|_X\le \Big|\int_{t'}^{t''}\|b(s,x(s))\|_Xds\Big|\le M|t'-t''|\to 0$$
as $t',t''\to t^*.$

Thus $x(t)$ satisfies the Cauchy criterion, and the limit (\ref{zcxv1}) exists.

It is important to stress that, by definition, the solution $x(t)$ remains in $Q$ for all $t\in[0,t^*)$;
therefore $x_-\in \overline Q$.

Suppose that $x_-\in Q$. Then let $y(t)$ be the a solution to  problem (\ref{5y}) with initial condition $y(t^*)=x_-$.
From the existence and uniqueness theorem, we know that
$y(t)$  is defined for $t\in[t^*,T')$  when $T'-t^*>0$ is sufficiently small.

One can easily verify that the function
$$\tilde x(t) =
\left\{
    \begin{array}{ll}
        x(t), & \text{if } t\in[0,t^*)\\
        y(t), & \text{if } t\in[t^*,T')
    \end{array}\right.
    $$
is a solution to (\ref{5y}) and this solution extends $x(t)$   to the right beyond $t^*$.

Indeed, $\tilde x$ satisfies the following integral equation

\begin{align*}
\tilde x(t)&=x_-+\int_{t^*}^tb(s,y(s))ds\\
&=x(0)+\int_{0}^{t^*}b(s,x(s))ds+\int_{t^*}^tb(s,y(s))ds\\
&=x(0)+\int_{0}^{t}b(s,\tilde x(s))ds,\quad t\ge t^*.\end{align*}

This proves  the theorem.

\begin{rem}\label{xdvfcrt00}
In the case where $t^*=T$  the limit (\ref{zcxv1})  exists as well,  and $x_-\in\overline Q.$\end{rem}


\begin{thebibliography}{99}


\bibitem{aga} R. P. Agarwal, M. Meehan, D. O'Regan:   Fixed Point Theory and Applications. Cambridge University Press (2001).
\bibitem{[1]} J. C. Alexander and J. A. Yorke: The homotopy continuation method: numerically implementable topological procedures. Transactions of the American Mathematical Society, vol. 242, pp. 271–284, 1978.

\bibitem{[2]} E. L. Allgower: A survey of homotopy methods for smooth mappings. in: E. L. Allgower, K. Glashoff, and H. O. Peitgen (eds.), Numerical Solution of Nonlinear Equations, Springer-Verlag, Berlin, pp. 2–29, 1981.
\bibitem{man} Alamri, B. Solving Integral
Equation and Homotopy Result via
Fixed Point Method. Mathematics
2023, 11, 4408.

\bibitem{[3]} E. L. Allgower and K. Georg, Numerical Continuation Methods: An Introduction, Springer Series in Computational Mathematics, vol. 13, Springer-Verlag, New York, 1990.

\bibitem{[4]} E. L. Allgower, K. Glashoff, and H. O. Peitgen (eds.), Proceedings of the Conference on Numerical Solution of Nonlinear Equations (Bremen, July 1980), Lecture Notes in Mathematics, vol. 878, Springer-Verlag, Berlin, 1981.
\bibitem{AR} A.V. Arutyunov and S.E. Zhukovskiy: Solvability of nonlinear degenerate equations and estimates for inverse functions. Sbornik: Mathematics, 216(1) (2025).


\bibitem{[5]} S. Bernstein, "Sur la g\'en\'eralisation du probl\'eme de Dirichlet," Mathematische Annalen, vol. 62, pp. 253–271, 1906.

\bibitem{brezis} H. Brezis: Functional Analysis, Sobolev Spaces and Partial Differential Equations. Springer Science \& Business Media, New York, 2011.
\bibitem{gran} A. Granas and J. Dugundji: Fixed Point Theory, Springer Monographs in Mathematics, Springer-Verlag, New York, 2003.

\bibitem{gar}
                        Sandoval-Hernandez, Mario A., Jimenez-Islas, Hugo, Vazquez-Leal, Hector, Quemada-Villagómez, Miriam L. and Lopez-Gonzalez, María de la Luz. "Exploring homotopy with hyperspherical tracking to find complex roots with application to electrical circuits" Open Mathematics, vol. 23, no. 1, 2025, pp. 20240115.

\bibitem{safonov}M. V. Safonov: The abstract Cauchy–Kovalevskaya theorem in a weighted Banach space. Communications on Pure and Applied Mathematics, 48(6), 629–637.

\bibitem{hart}P. Hartman: Ordinary Differential Equations. Birkh\"auser, Boston, 1982.

\bibitem{mil}V. M. Millionschikov: On the theory of differential equations in locally convex spaces. Matematicheskii Sbornik, 51(93):1, 105–114  (1960).
\bibitem{nir} L. Nirenberg: Topics in Nonlinear Functional Analysis. New York, 1974.
\bibitem{re}W. C. Rheinboldt: Local mapping relations and global implicit function theorems. Transactions of the American Mathematical Society, vol. 138, pp. 183–198, 1969.
\bibitem{gar}
                        Sandoval-Hernandez, Mario A., Jimenez-Islas, Hugo, Vazquez-Leal, Hector, Quemada-Villag\'omez, Miriam L. and Lopez-Gonzalez, Maria de la Luz. "Exploring homotopy with hyperspherical tracking to find complex roots with application to electrical circuits" Open Mathematics, vol. 23, no. 1, 2025.pp. 20240115.
\bibitem{swarz} L. Schwartz: Analyse Math\'ematique. Hermann, 1967.
\bibitem{ang}J. M. Soriano and V. G. Angelov: A Zero of a Proper Mapping. Fixed Point Theory, Vol. 4, No. 1, 2003, pp. 97–103.
\bibitem{[22]} J. M. Soriano: Zeros of Compact Perturbations of Proper Mappings. Applied Nonlinear Analysis, vol. 7, no. 4, pp. 31–37, 2000.

\bibitem{[23]} J. M. Soriano: Fredholm and Compact Mappings Sharing a Value. Applied Mathematics and Mechanics (English Edition), vol. 24, no. 8, pp. 941–948, 2003.


\bibitem{[24]} J. M. Soriano: Compact Mapping and Proper Mapping Between Banach Spaces that Share a Value. Mathematica Balkanica (New Series), vol. 14, no. 1-2, pp. 161–166, 2000.

\bibitem{Yak}M. N. Yakovlev: The solutions of systems of non-linear
equations by a method of differentiation with respect to a
parameter,
Zh. Vychisl. Mat. Mat. Fiz., 1964, Volume 4, Number 1, 146–149.





\bibitem{zy} O. Zubelevich: Abstract version of the Cauchy–Kowalevski Problem CEJM 2(3) 2004, p. 382–387.

 \end{thebibliography}
\end{document}